\font\script=eusm10.
\font\sets=msbm10.
\font\stampatello=cmcsc10.
\font\symbols=msam10.

\def\sgn{\hbox{\rm sgn}}

\def\defineq{\buildrel{def}\over{=}}
\def\defin{\buildrel{def}\over{\Longleftrightarrow}}
\def\square{\hbox{\vrule\vbox{\hrule\phantom{s}\hrule}\vrule}}
\def\1{\hbox{\bf 1}}
\def\C{\hbox{\sets C}}
\def\N{\hbox{\sets N}}

\def\R{\hbox{\sets R}}
\def\Z{\hbox{\sets Z}}
\def\doublesum{\mathop{\sum \sum}}

\def\integrale{\mathop{\int}}

\def\Corr{\hbox{\script C}}

\def\EssBdd{\hbox{\symbols n}_{\varepsilon}\,}

\par
\centerline{\bf ON THE CORRELATIONS, SELBERG INTEGRAL AND SYMMETRY}
\centerline{\bf OF SIEVE FUNCTIONS IN SHORT INTERVALS}
\bigskip
\centerline{by G.Coppola}
\bigskip
{
\font\eightrm=cmr8
\eightrm {
\par
{\bf Abstract.} We study the arithmetic (real) function $f=g\ast \1$, with $g$ \lq \lq essentially bounded\rq \rq \thinspace and supported over the integers of $[1,Q]$. In particular, we obtain non-trivial bounds, through $f$ \lq \lq correlations\rq \rq, for the \lq \lq Selberg integral\rq \rq \thinspace and the \lq \lq symmetry integral\rq \rq \thinspace of $f$ in almost all short intervals $[x-h,x+h]$, $N\le x\le 2N$, beyond the \lq \lq classical\rq \rq \thinspace level, up to level of distribution, say, $\lambda=\log Q/\log N < 2/3$ (for enough large $h$). This time we don't apply Large Sieve inequality, as in our paper [C-S]. Precisely, our method is completely elementary. 
}
\footnote{}{\par \noindent {\it Mathematics Subject Classification} $(2000) : 11N37, 11N25.$}
}
\bigskip
\par
\noindent {\bf 1. Introduction and statement of the results.}
\smallskip
\par
We study \lq \lq {\stampatello sieve functions}\rq \rq \thinspace, i.e. real arithmetic functions $f=g\ast \1$ (see hypotheses on $g$ in the sequel), in almost all the short intervals $[x-h,x+h]$ (i.e., almost all stands $\forall x\in [N,2N]$, except $o(N)$ of them and short means, say, $h\to \infty$ and $h=o(N)$, as $N\to \infty$). Here, as usual, $\1(n)=1$ is the constant-$1$ arithmetic function and $\ast$ is the Dirichlet product (esp., [T]). In order to study the sum of $f$ values in a.a. (abbreviates almost all, now on) the intervals $[x-h,x+h]$, we define (in analogy with the classical Selberg integral, see [C-S]) the \lq \lq {\stampatello Selberg integral}\rq \rq \thinspace of \thinspace $f$ \thinspace as:\enspace $J_f(N,h)\defineq {\displaystyle \int_{N}^{2N} \Big| \sum_{0<|n-x|\le h} } f(n) - M_f(2h)\Big|^2\,dx$, \thinspace where (from heuristics in accordance with the classical case) we expect the \lq \lq mean-value\rq \rq \thinspace to be \thinspace $M_f(2h)\defineq 2h\sum_{d}g(d)/d$ 
\par
\noindent
(that converges in interesting cases and under our hypotheses on $g$, see the sequel; also, $d\le 2N+h$, here).
\par
\noindent
Furthermore, this definition comes from what the \lq \lq natural\rq \rq \thinspace choice of $M_f(2h)$ \thinspace is (recall $[\enspace]=$ {\stampatello integer part}): 
$$
2h\Big({1\over x}\sum_{n\le x}f(n)\Big) = {{2h}\over x}\sum_{d}g(d)\left[ {x\over d}\right] = 2h\sum_{d}{{g(d)}\over d} + {\cal O}\Big( {h\over x}\sum_{d\le Q}|g(d)|\Big), 
$$
\par
\noindent
in fact, when $f=g\ast \1$, $g(q)=0$ for $q>Q$. Assuming $Q$ smaller than $x$ (in the sequel), we recover \thinspace $M_f(2h)$. 
\par
\noindent
Selberg integral counts the values of $f$ in a.a. $[x-h,x+h]$. We study their symmetry through the \lq \lq {\stampatello symmetry integral}\rq \rq \thinspace of $f$ (here $\sgn(0)\defineq 0$, $\sgn(r)\defineq {{|r|}\over r}$, $\forall r\neq 0$):\enspace $I_f(N,h)\defineq {\displaystyle \int_{N}^{2N} \Big| \sum_{|n-x|\le h} }\sgn(n-x)f(n)\Big|^2\,dx$. 
\par
We'll generalize the results given in [C-S] for these integrals, applying the Large Sieve inequality, in the case $g=\1$ of the divisor function $d=\1 \ast \1$. We point out that the procedure given there works, as well, for more general $g$ to bound $I_f$; but fails in the case of $J_f$, whenever $g$ is not constant (i.e., the Dirichlet \lq \lq flipping\rq \rq \thinspace of the divisors can't be applied). Here, we give another approach valid for both integrals, even for non-constant $g$. It is based on the \lq \lq correlations\rq \rq \thinspace of $f$. The {\stampatello correlation of} $f$ {\stampatello is} defined as ($\forall a\in \Z, a\neq 0$)
$$
\Corr_f(a)\defineq \sum_{n\sim N}f(n)f(n-a) 
= \sum_{\ell |a}\doublesum_{(d,q)=1}g(\ell d)g(\ell q){1\over q}\left( \left[ {{2N}\over {\ell d}}\right] - \left[ {N\over {\ell d}}\right]\right) + R_f(a) 
$$
\par
\noindent
(hereon \thinspace $x\sim X$ \thinspace is \thinspace $X<x\le 2X$), where, through the orthogonality of additive characters [V] as in Lemma 3  
\par
\noindent
(as usual, we will always write \enspace $e(\theta)\defineq e^{2\pi i \theta}$, $\forall \theta \in \R$ \thinspace and \thinspace $e_q(m)\defineq e(m/q)$, $\forall q\in \N$, $\forall m\in \Z$), say, 
$$
R_f(a)\defineq \sum_{\ell |a}\doublesum_{(d,q)=1}g(\ell d)g(\ell q){1\over q}\sum_{j\neq 0}e_q(-ja/\ell)\sum_{m\sim {N\over {\ell d}}}e_q(jdm) 
$$
\par
\noindent
(here, and in the following, \thinspace $j\neq 0$ \thinspace means that $j$ describes exactly once all classes $(\!\bmod \,\, q)$, except \thinspace $j\equiv 0(q)$); 
\par
\noindent
and $I_f$ is a sum (see Lemma 1) of these correlations, weighted with $W$ (name from the shape), $W$ {\stampatello even}, 
$$
W(a)\defineq \cases{2h-3a & if $0\le a\le h$\cr a-2h & if $h\le a\le 2h$\cr 0 & if $a>2h$\cr} 
\quad \enspace \Longrightarrow \enspace \quad \sum_{a\in \Z}W(a)=0.
$$
\par
\noindent
In complete analogy, Lemma 2 gives the Selberg integral $J_f(N,h)$ as a weighted sum of correlations, with (Selberg) weight \thinspace $S(a)\defineq \max(2h-|a|,0)$. Notice that $S$ is always non-negative (while $W$ oscillates in sign).

\vfil
\eject
\par
\noindent
(Here, as usual, $F=o(G)\defin \lim {F/G}=0$ and $F={\cal O}(G) \defin \exists c>0: |F|\le c\,G$ are Landau's notation. Also, when $c$ depends on $\varepsilon$, we'll write $F={\cal O}_{\varepsilon}(G)$ or, like Vinogradov, $F\ll_{\varepsilon}G$). We call an arithmetical function {\stampatello essentially bounded} when, $\forall \varepsilon>0$, its $n-$th value is at most ${\cal O}_{\varepsilon}(n^{\varepsilon})$ and we'll write $\EssBdd 1$; i.e.,
$$
F(N)\EssBdd G(N) \defin \forall \varepsilon>0 \enspace F(N)\ll_{\varepsilon} N^{\varepsilon} G(N)\quad (\hbox{\rm as} \; N\to \infty)
$$
\par
\noindent
e.g., the divisor function $d(n)$ is essentially bounded (like many other number-theoretic $f$) and we remark that $f=g\ast \1$ is essentially bounded if and only if $g$ is (from M\"{o}bius inversion, see [D]). From Lemma 2, applying Lemma 3 to $f$ correlations, together with \thinspace $\sum_a S(a\ell)=4h^2/\ell + {\cal O}(h)$, uniformly $\forall \ell \in \N$ (like in (1), see Lemma 4 proof), we get 
$$
J_f(N,h) = \sum_{\ell \le 2h}\doublesum_{(d,q)=1}g(\ell d)g(\ell q)\sum_a S(a\ell)R_f(a) + {\cal O}_{\varepsilon}\left(N^{\varepsilon}\left(Nh+Qh^2\right)\right) 
$$
\par
\noindent
In fact, (compare the discussion about $M_f(2h)$, above)
$$
f\EssBdd 1 \Rightarrow M_f(2h)\EssBdd h, \enspace \sum_{n\sim N}f(n)=\sum_{d}g(d)\left(\left[ {{2N}\over d}\right] - \left[ {N\over d}\right]\right) = N\sum_{d}{{g(d)}\over d} + {\cal O}_{\varepsilon}\left( N^{\varepsilon}Q\right). 
$$
\par
\noindent
We recall\enspace $\Vert r\Vert \defineq {\displaystyle \min_{n\in \Z} }|r-n|$ is the {\stampatello distance from integers}. We abbreviate \thinspace $n\equiv a(\!\bmod\, q)$ \thinspace with \thinspace $n\equiv a(q)$. 
\medskip
\par
We give our main result. 
\smallskip
\par
\noindent {\stampatello Theorem.} {\it Let } $N,h,Q\in \N$, {\it be such that } $h\to \infty$, $Q\ll N$ {\it and } $h=o(N)$, {\it as } $N\to \infty$. {\it Let } $f:\N \rightarrow \R$ {\it be essentially bounded, with } $f=g \ast \1$ {\it and } $g(q)=0$ $\forall q>Q$. {\it Then}
$$
J_f(N,h)\EssBdd Nh + h^3 + Q^2 h + Qh^2; \qquad I_f(N,h)\EssBdd Nh + h^3 + Q^2 h + Qh^2.
$$
\par
\noindent
{\it Also, only for the symmetry integral } $I_f(N,h)$, 
$$
I_f(N,h) = 2\sum_a S(a)\left( \Corr_f(a) - \Corr_f(a+h)\right) + {\cal O}_{\varepsilon}\left( N^{\varepsilon}(Nh+h^3) \right). 
$$
\smallskip
\par
\noindent
{\bf Remark.} We explicitly point out that our Theorem implies non-trivial estimates \thinspace $J_f(N,h)\ll {{Nh^2}\over {N^{\varepsilon}}}$ and $I_f(N,h)\ll {{Nh^2}\over {N^{\varepsilon}}}$ \thinspace for both integrals, with {\stampatello level of distribution}, say, \thinspace ${{\log Q}\over {\log N}} \defineq \lambda < {{1+\theta}\over 2}$, where, say, \thinspace $\theta \defineq (\log h)/(\log N)$ is the {\stampatello width}; hence, level up to $2/3$, when the width is above $1/3$. (The same  result can also be achieved with the method of [C-S], but only for $I_f$.) 
\medskip
\par
In fact, an immediate consequence of our Theorem is the following 
\smallskip
\par
\noindent {\stampatello Corollary.} {\it Let } $0<\theta<1$, $0\le \lambda<{{1+\theta}\over 2}$ {\it and } $N,h,Q\in \N$, {\it be such that } $N^{\theta}\ll h\ll N^{\theta}$, $N^{\lambda}\ll Q\ll N^{\lambda}$, {\it as } $N\to \infty$. {\it Let } $f:\N \rightarrow \R$ {\it be essentially bounded, with } $f=g \ast \1$ {\it and } $g(q)=0$ $\forall q>Q$. {\it Then} $\exists \varepsilon_0=\varepsilon_0(\theta,\lambda)>0$ ({\it depending only on } $\theta,\lambda$) {\it such that}
$$
J_f(N,h)\ll_{\varepsilon_0} Nh^2N^{-\varepsilon_0}, \qquad I_f(N,h)\ll_{\varepsilon_0} Nh^2N^{-\varepsilon_0}.
$$
\bigskip
\par
\noindent
The paper is organized as follows: 
\medskip
\par
\noindent $\diamond$ we will give our Lemmas in the next section; 
\smallskip
\par
\noindent $\diamond$ then we will prove our Theorem in section 3. 

\vfil
\eject
\par
\noindent {\bf 2. Lemmas.}
\smallskip
\par
\noindent {\bf Lemma 1.} {\it Let } $N,h\in \N$, {\it with } $h\to \infty$ {\it and } $h=o(N)$ {\it as } $N\to \infty$. 
{\it If } $f:\N \rightarrow \R$ {\it has } $\Vert f\Vert_{\infty}={\displaystyle \max_{n\le 2N+h} }|f(n)|$, 
$$
\int_{N}^{2N} \Big| \sum_{|n-x|\le h}\sgn(n-x)f(n)\Big|^{2}\,dx 
= \sum_{a}W(a)\Corr_f(a) + {\cal O}\left( h^3\Vert f\Vert_{\infty}^2 \right). 
$$
\smallskip
\par
\noindent {\it Proof.}$\!$ This is a kind of dispersion method, without \lq \lq expected mean\rq \rq :the main term \lq \lq vanishes\rq \rq. Use $f$ {\stampatello real}:
$$
\qquad \enspace 
I_f(N,h)=D_f(N,h)+2\doublesum_{N-h < n_1 < n_2 \le 2N+h}f(n_1)f(n_2)\integrale_{x\sim N,|x-n_1|\le h,|x-n_2|\le h}\sgn(x-n_1)\sgn(x-n_2)dx;
$$
\par
here ($I_f$ is the integral above and) \enspace $D_f(N,h):= {\displaystyle \sum_{N-h<n\le 2N+h}f^2(n)\integrale_{N<x\le 2N,0<|x-n|\le h}dx } =$
$$
\qquad 
= \sum_{N+h<n\le 2N-h}f^2(n)\integrale_{|x-n|\le h}dx +{\cal O}\left(h\Vert f\Vert_{\infty}^2 \left(\sum_{|n-N|\le h}1+\sum_{|n-2N|\le h}1\right)\right) 
= W(0)\Corr_f(0)+{\cal O}\left(h^2 \Vert f\Vert_{\infty}^2 \right) 
$$
\par
is the {\stampatello diagonal}. The remainder, here, is (negligible) in the second one. Since (for $a>0$) 
\par
$$
\Corr_f(-a) = \sum_{n\sim N}f(n)f(n+a)= \sum_{N+a<m\le 2N+a}f(m-a)f(m) = \Corr_f(a)+{\cal O}\left(a\Vert f\Vert_{\infty}^2 \right),
$$ 
\par
$W$ {\stampatello even} and \thinspace $W(a)\ll h$ \thinspace $\Rightarrow$ \thinspace ${\displaystyle \sum_{0<a\le 2h}W(a)\Corr_f(-a)=\sum_{0<a\le 2h} }W(a)\Corr_f(a)+{\cal O}\left(h^3 \Vert f\Vert_{\infty}^2 \right)$, we confine to: 
$$
\enspace 
I_f(N,h)-D_f(N,h):=W(0)\Corr_f(0)+2\sum_{0<a\le 2h}W(a)\Corr_f(a)+E_f(N,h), \thinspace \hbox{\rm say,} \thinspace E_f(N,h)\ll h^3\Vert f\Vert_{\infty}^2.
\leqno{(*)}
$$
\par
The left-hand side, changing variables, namely $n=n_1$, $a=n_2-n_1$, $s=x-n_1$, is (introducing the 
\par
remainders which shall take part of the final \thinspace $E_f(N,h)$, here)
\vskip-0.15cm
$$
2\sum_{N-h<n<2N+h}f(n)\sum_{0<a\le 2h,a\le 2N+h-n}f(n+a)\integrale_{N-n<s\le 2N-n,|s|\le h,|s-a|\le h}\sgn(s)\sgn(s-a)ds =
$$
$$
= 2\sum_{N-h<n\le 2N-h}f(n)\sum_{0<a\le 2h}f(n+a)\integrale_{s>N-n,|s|\le h,|s-a|\le h}\sgn(s)\sgn(s-a)ds+E_1 = 
$$
$$
= 2\sum_{N+h<n\le 2N-h}f(n)\sum_{0<a\le 2h}f(n+a)\integrale_{|s|\le h,|s-a|\le h}\sgn(s)\sgn(s-a)ds+E_1+E_2 = 
$$
$$
= 2\sum_{N<n\le 2N}f(n)\sum_{0<a\le 2h}f(n+a)W(a)+E_1+E_2+E_3,
\enspace W(a)\!:=\!\!\!\integrale_{{|s|\le h}\atop {|s-a|\le h}}\!\!\sgn(s)\sgn(s-a)ds\!\ll h, 
$$
\par
whence\enspace $E_3\ll\!\!{\displaystyle \left(\sum_{N<n\le N+h}+\sum_{2N-h<n\le 2N}\right)|f(n)|\sum_{0<a\le 2h} }\!|f(n+a)|h\ll h^3 \Vert f\Vert_{\infty}^2$ \enspace is a \thinspace \lq \lq {\stampatello tail}\rq \rq, like: 
$$
E_1\ll \sum_{|n-2N|\le h}|f(n)|\sum_{0<a\le 2h}|f(n+a)|h, \thinspace 
E_2\ll \sum_{|n-N|\le h}|f(n)|\sum_{0<a\le 2h}|f(n+a)|h \enspace \hbox{\rm are} \thinspace \ll h^3 \Vert f\Vert_{\infty}^2.\enspace \square
$$

\vfil
\eject

\par
\noindent {\bf Lemma 2.} {\it Let } $N,h\in \N$, {\it with } $h\to \infty$ {\it and } $h=o(N)$ {\it as } $N\to \infty$. 
{\it If } $f:\N \rightarrow \R$ {\it has } $\Vert f\Vert_{\infty}={\displaystyle \max_{n\le 2N+h} }|f(n)|$, 
$$
\int_{N}^{2N} \Big| \sum_{0<|n-x|\le h}f(n)-M_f(2h)\Big|^{2}\,dx  
= \sum_a S(a)\Corr_f(a) - 4M_f(2h)h\sum_{n\sim N}f(n) + M_f^2(2h)N 
+ 
$$
$$
+ {\cal O}\left( h^3\, \Vert f\Vert_{\infty}^2 + h^2 \Vert f\Vert_{\infty} |M_f(2h)|\right). 
$$
\smallskip
\par
\noindent {\it Proof.}$\!$ This is a direct application of dispersion method [L]. Use $f$ {\stampatello real} (ignoring, now, sets of measure zero):
$$
\quad \enspace 
J_f(N,h)=D_f(N,h)+2\doublesum_{N-h < n_1 < n_2 \le 2N+h}f(n_1)f(n_2)\integrale_{x\sim N,|x-n_1|\le h,|x-n_2|\le h}dx -
$$
$$
- 2 M_f(2h)\sum_{N-h<n\le 2N+h}f(n)\integrale_{x\sim N,|x-n|\le h}dx + M_f^2(2h)\int_{N}^{2N}dx = 
$$
$$
\qquad 
= D_f(N,h)+2\doublesum_{N-h < n_1 < n_2 \le 2N+h}f(n_1)f(n_2)\integrale_{x\sim N,|x-n_1|\le h,|x-n_2|\le h}dx 
- 4h M_f(2h)\sum_{n\sim N}f(n) + M_f^2(2h)N, 
$$
\par
save an error which is \enspace ${\cal O}(|M_f(2h)|h^2\Vert f\Vert_{\infty})$; here ($J_f$ is the integral above and) 
$$
D_f(N,h)= \sum_{N-h<n\le 2N+h}f^2(n)\integrale_{x\sim N,0<|x-n|\le h}dx = 
$$
$$
\qquad 
= \sum_{N<n\le 2N}f^2(n)\integrale_{0<|x-n|\le h}dx +{\cal O}\left(h\Vert f\Vert_{\infty}^2 \left(\sum_{|n-N|\le h}1+\sum_{|n-2N|\le h}1\right)\right) 
= S(0)\Corr_f(0)+{\cal O}\left(h^2 \Vert f\Vert_{\infty}^2 \right) 
$$
\par
is the same diagonal (with same negligible remainder) of Lemma 1. In fact, we closely follow its proof; 
\par
due to: $S$ {\stampatello even} and \thinspace $S(a)\ll h$ \thinspace $\Rightarrow$ \thinspace ${\displaystyle \sum_{0<a\le 2h}S(a)\Corr_f(-a)=\sum_{0<a\le 2h} }S(a)\Corr_f(a)+{\cal O}\left(h^3 \Vert f\Vert_{\infty}^2 \right)$, we confine to 
$$
\doublesum_{N-h < n_1 < n_2 \le 2N+h}f(n_1)f(n_2)\integrale_{x\sim N,|x-n_1|\le h,|x-n_2|\le h}dx 
- \sum_{0<a\le 2h}S(a)\Corr_f(a) := E_f(N,h)
\buildrel{\hbox{\rm say}}\over{\ll}
h^3\Vert f\Vert_{\infty}^2.
\leqno{(*)}
$$
\par
The left-hand side, changing variables, namely $n=n_1$, $a=n_2-n_1$, $s=x-n_1$, is (see Lemma 1 proof) 
$$
\sum_{N-h<n<2N+h}f(n)\sum_{0<a\le 2h,a\le 2N+h-n}f(n+a)\integrale_{N-n<s\le 2N-n,|s|\le h,|s-a|\le h}ds =
$$
$$
= \sum_{N<n\le 2N}f(n)\sum_{0<a\le 2h}f(n+a)S(a)+E_1+E_2+E_3,
\enspace S(a)\!:=\!\!\!\integrale_{{|s|\le h}\atop {|s-a|\le h}}ds\ll h, 
$$
\par
whence\enspace $E_3\ll\!\!{\displaystyle \left(\sum_{N<n\le N+h}+\sum_{2N-h<n\le 2N}\right)|f(n)|\sum_{0<a\le 2h} }\!|f(n+a)|h\ll h^3 \Vert f\Vert_{\infty}^2$ \enspace is a \thinspace \lq \lq {\stampatello tail}\rq \rq, like: 
$$
\enspace \enspace 
E_1\ll \sum_{|n-2N|\le h}|f(n)|\sum_{0<a\le 2h}|f(n+a)|h, \enspace \thinspace 
E_2\ll \sum_{|n-N|\le h}|f(n)|\sum_{0<a\le 2h}|f(n+a)|h \enspace \hbox{\rm are} \thinspace \ll h^3 \Vert f\Vert_{\infty}^2.\enspace \square
$$
\par
\noindent {\bf Lemma 3.} {\it Let } $N,h,Q\in \N$, {\it where } $h\to \infty$, $h=o(N)$ \thinspace {\it and } $Q\ll N$, {\it as } $N\to \infty$. {\it Let } $f = g\ast \1$, {\it where } $g:\N \rightarrow \R$, {\it with } $q>Q \enspace \Rightarrow \enspace g(q)=0$. {\it Then} 
$$
a\neq 0 \Rightarrow \Corr_f(a) = \sum_{\ell | a}\doublesum_{(d,q)=1}g(\ell d)g(\ell q){1\over q}\left( \left[{{2N}\over {\ell d}}\right] - \left[{N\over {\ell d}}\right]\right)+ R_f(a), \enspace \hbox{\it where, say, as in the introduction}
$$
$$
R_f(a) = \sum_{\ell | a}\doublesum_{(d,q)=1}g(\ell d)g(\ell q){1\over q}\sum_{j\neq 0}e_q(-ja/\ell)\sum_{m\sim {N\over {\ell d}}}e_q(jdm), \qquad \forall a\neq 0.
$$
\par
\noindent
{\it Also, every weight function } $K:\N \rightarrow \C$, $K$ {\stampatello even}, {\it with } $K(0)=2h$, {\it gives} 
$$
\sum_{a}K(a)R_f(a) 
= \sum_{\ell \le 2h}\doublesum_{(d,q)=1}g(\ell d)g(\ell q){1\over q}\sum_{j\neq 0}\sum_{m\sim {N\over {\ell d}}}\cos {{2\pi jdm}\over q}\sum_{a\neq 0}K(a\ell)e_q(ja) + 2h\Corr_f(0).
$$
\par
\noindent {\it Proof.} We'll always assume $a$ non-zero. First of all, we start from the correlation, that is: 
$$
\Corr_f(a) = \sum_{n\sim N}f(n)f(n-a) 
= \sum_{d}\sum_{q}g(d)g(q)\sum_{{N<n\le 2N}\atop {{n\equiv 0(d)}\atop {n\equiv a(q)}}}1 
= \doublesum_{(d,q)|a}g(d)g(q)\sum_{{{N\over d}<m\le {{2N}\over d}}\atop {md\equiv a (q)}}1, 
$$
\par
since last congruence is solveable if and only if the GCD $(d,q)$ divides $a$; \enspace changing variables, this is 
$$
\sum_{\ell | a}\doublesum_{(d,q)=1}g(\ell d)g(\ell q)\sum_{{{N\over {\ell d}}<m\le {{2N}\over {\ell d}}}\atop {md\equiv {a\over {\ell}} (q)}}1 
= \sum_{\ell | a}\doublesum_{(d,q)=1}g(\ell d)g(\ell q){1\over q}\left( \left[ {{2N}\over {\ell d}}\right] - \left[ {N\over {\ell d}}\right]\right) + 
R_f(a), 
$$
\par
using the orthogonality of additive characters (see [V]): here \thinspace $R_f(a)$ \thinspace is as above; summing on $a$ with $K$, 
$$
\sum_{a\neq 0}K(a)R_f(a) 
= \sum_{\ell \le 2h}\doublesum_{(d,q)=1}g(\ell d)g(\ell q){1\over q}\sum_{j\neq 0}\sum_{m\sim {N\over {\ell d}}}e_q(jdm)\sum_{b\neq 0}K(b\ell)e_q(jb) = 
$$
\par
(using $K$ even, here)
$$
= \sum_{\ell \le 2h}\doublesum_{(d,q)=1}g(\ell d)g(\ell q){1\over q}\sum_{j\neq 0}\sum_{m\sim {N\over {\ell d}}}e_q(jdm)\sum_{a\neq 0}K(a\ell)\cos{{2\pi ja}\over q} = 
$$
$$
= \sum_{\ell \le 2h}\doublesum_{(d,q)=1}g(\ell d)g(\ell q){1\over q}\sum_{j\neq 0}\sum_{m\sim {N\over {\ell d}}}\cos{{2\pi jdm}\over q}\sum_{a\neq 0}K(a\ell)\cos{{2\pi ja}\over q} = 
$$
$$
= \sum_{\ell \le 2h}\doublesum_{(d,q)=1}g(\ell d)g(\ell q){1\over q}\sum_{j\neq 0}\sum_{m\sim {N\over {\ell d}}}\cos{{2\pi jdm}\over q}\sum_{a\neq 0}K(a\ell)e_q(ja).
$$
\par
(We used once more $K$ even, here.) Then, the thesis, adding the term \thinspace $K(0)\Corr_f(0)=2h\Corr_f(0).\enspace \square$
\bigskip
\par
\noindent {\bf Remark.} We explicitly point out that, in our hypotheses on $f$ (i.e., real and essentially bounded) 
$$
2h\Corr_f(0) = 2h\sum_{n\sim N}f^2(n) 
\EssBdd Nh, 
$$
\par
\noindent
a trivial estimate which will be useful in future occurrences. 

\vfil
\eject

\par
\noindent {\bf Lemma 4.} {\it Defining}, $\forall h\in \N$, {\it the weight } $W$ {\it as above, we have}, $\forall q\in \N$, $\forall \beta \not \in \Z$, $\forall \ell \in \N$, $\forall \alpha \in \R$,
$$
\sum_{{a}\atop {a\equiv 0(\!\!\bmod q)}}W(a)=2q\left\Vert {h\over q}\right\Vert;
\quad
\sum_{0\le |a|\le 2h}W(a)e(a\beta) = {{4\sin^4(\pi h \beta)}\over {\sin^2(\pi \beta)}};
\quad
\sum_{b}W(\ell b)e(b\alpha)\ge 0.
\leqno{(1)}
$$
\par
\noindent
{\it Also, more in general} ({\it in the same hypotheses}), {\it abbreviating } $E_{X}(\beta)\defineq \sum_{0\le |a|\le X}e(a\beta)$, {\it we have} 
$$
{1\over {\ell}}\sum_{a}W(a\ell)e(a\beta) = {{4\sin^2 \pi \beta [{h\over {\ell}}]-\sin^2 \pi \beta [{{2h}\over {\ell}}]}\over {\sin^2 \pi \beta}} + 
4\left\{ {h\over {\ell}}\right\} E_{h\over {\ell}}(\beta) - \left\{ {{2h}\over {\ell}}\right\} E_{{2h}\over {\ell}}(\beta) = 
\leqno{(2)}
$$
$$
= 2\left(1-\cos \left(2\pi \beta \left[{h\over {\ell}}\right]\right)\right)\sum_{0\le |a|\le {h\over {\ell}}}\left( \left[{h\over {\ell}}\right]-|a|\right) e(a\beta) - \left( 2\left\{ {h\over {\ell}}\right\} - \left\{ {{2h}\over {\ell}}\right\} \right)E_{2[{h\over {\ell}}]}(\beta) + 
$$
$$
+ 4\left\{ {h\over {\ell}}\right\} E_{h\over {\ell}}(\beta) - \left\{ {{2h}\over {\ell}}\right\} E_{{2h}\over {\ell}}(\beta).
$$
\smallskip
\par
\noindent {\it Proof.}$\!$ Hereon $n\le X$ in a sum means $1\le n\le X$. We will prove (1), even if it's a special case of (2); 
\vskip-0.25cm
$$
\sum_{{a}\atop {a\equiv 0(\!\!\bmod q)}}W(a) 
= 2h+4h\left( \left[ {h\over q}\right]-\left[ {{2h}\over q}\right]+\left[ {h\over q}\right]\right)+2q\left( \sum_{{h\over q}<b\le {{2h}\over q}}b - 3\sum_{b\le {h\over q}}b\right)
$$
\vskip-0.25cm
$$
\!=q\left( \left\{ {{2h}\over q}\right\}^2 - 4\left\{ {h\over q}\right\}^2 - \left\{ {{2h}\over q}\right\} + 4\left\{ {h\over q}\right\} \right).
$$
\par
Using \thinspace $\forall \alpha\in \R$ \thinspace that \thinspace $\{ 2\alpha \} = \{ 2\{ \alpha \} \}=\cases{2\{ \alpha \} & if \thinspace $0\le \{ \alpha \}<1/2$\cr 2\{ \alpha \}-1 & if \thinspace $1/2 \le \{ \alpha \}<1$\cr}$ \enspace we get the first.
\par
We come, now, to the second: \enspace ${\displaystyle \sum_{0\le |a|\le 2h}W(a)e(a\beta)=2h+2\sum_{a\le 2h} }W(a)\cos {2\pi a\beta}=2h+2\Sigma$, \enspace say; then, 
\par
partial summation gives \thinspace ${\displaystyle \Sigma=4\sum_{a\le h} C_a(\beta) - 4C_h(\beta) - \sum_{a\le 2h} }C_a(\beta) + C_{2h}(\beta)$, \thinspace say, where\enspace $\forall X\in \N$, $\forall \theta \not \in \Z$ \enspace
\vskip-0.15cm
$$
C_X(\theta)\defineq \sum_{n\le X}\cos(2\pi n\theta)={{\sin(2\pi \theta X)}\over {2\tan(\pi \theta)}}-{{1-\cos(2\pi \theta X)}\over 2} \enspace (\hbox{a \thinspace well-known\thinspace formula}) 
$$
\vskip-0.15cm
\par
to get \enspace $\Sigma = {\displaystyle 2\cot(\pi \beta) \sum_{a\le h}\sin(2\pi a\beta) - 2h - 2C_h(\beta) - {1\over 2} \cot(\pi \beta)\sum_{a\le 2h}\sin(2\pi a\beta) + {1\over 2}C_{2h}(\beta) + h =}$
\vskip-0.10cm
$$
\,= \cot(\pi \beta) \left( 2S_h(\beta)-{1\over 2}S_{2h}(\beta)\right) - h - 2C_h(\beta) + {1\over 2}C_{2h}(\beta),\enspace \hbox{\rm say,}\; \forall X\in \N, \forall \theta \not \in \Z 
$$
\vskip-0.10cm
$$
S_X(\theta)\defineq \sum_{n\le X}\sin(2\pi n\theta)={{\sin^2(\pi \theta X)}\over {\tan(\pi \theta)}}+{{\sin(2\pi \theta X)}\over 2}
$$
\par
Then, since \enspace ${{1-\cos(2\pi \beta X)}\over 2}=\sin^2(\pi \beta X)$, both for \thinspace $X=h$ \thinspace and \thinspace $X=2h$,
$$
\Sigma = \cot^2(\pi \beta) \left(2\sin^2(\pi \beta h)-{1\over 2}\sin^2(2\pi \beta h)\right) + 2\sin^2(\pi \beta h) - {1\over 2}\sin^2(2\pi \beta h) - h = 
$$
$$
= 2\cot^2(\pi \beta)\left( 1 - \cos^2(\pi h\beta) \right)\sin^2(\pi h\beta) + 2\sin^2(\pi h\beta)\left( 1 - \cos^2(\pi h\beta) \right) - h = 
$$
$$
= 2 \left( \cot^2(\pi \beta) + 1 \right) \sin^4(\pi h\beta) - h 
= {{2\sin^4(\pi h\beta)}\over {\sin^2(\pi \beta)}} - h.
$$
\par
This gives the second. Finally, the third follows from: \thinspace $\forall \ell \in \N$
$$
\sum_b W(\ell b)e(b\alpha)\ge 0 \enspace \forall \alpha \in \R \Longleftrightarrow \sum_{a\equiv 0(\!\!\bmod \ell)} W(a)e(a\beta)\ge 0 \enspace \forall \beta \in \R
$$
\par
which, using the orthogonality of additive characters [V] and \enspace $\sum_a W(a)e(a\beta)\ge 0$ \thinspace $\forall \beta \in \R$, is 
$$
\sum_{a\equiv 0(\!\!\bmod \ell)} W(a)e(a\beta)={1\over {\ell}} \sum_{j\le \ell}\sum_a W(a)e(a\beta)e_{\ell}(ja)={1\over {\ell}} \sum_{j\le \ell}\sum_a W(a)e\left( a\left( \beta + {j\over {\ell}}\right) \right)\ge 0.
$$
\par
(We explicitly remark that this last property isn't \lq \lq visible\rq \rq \thinspace from (2): not an immediate consequence.) 
\medskip
\par
We come, now, to (2): \enspace ${\displaystyle \sum_{0\le |a|\le 2h}W(a\ell)e(a\beta)=2h+2\sum_{a\le {{2h}\over {\ell}}} }W(a\ell)\cos {2\pi a\beta}=2h+2\Sigma_{\ell}$, \enspace say; then, 
$$
{1\over {\ell}}\Sigma_{\ell}=4\sum_{a\le \left[{h\over {\ell}}\right]} C_a(\beta) - 4C_{\left[{h\over {\ell}}\right]}(\beta) - \sum_{a\le \left[{{2h}\over {\ell}}\right]}C_a(\beta) + C_{\left[{{2h}\over {\ell}}\right]}(\beta) + \left( 4\left\{ {h\over {\ell}}\right\}C_{\left[{h\over {\ell}}\right]}(\beta) - \left\{ {{2h}\over {\ell}}\right\}C_{\left[{{2h}\over {\ell}}\right]}(\beta)\right),
$$
\par
from partial summation (the term in brackets isn't present whenever $\ell=1$); then, (see above formulas)  
$$
{1\over {\ell}}\Sigma_{\ell}=\cot(\pi \beta)\left( 2S_{\left[{h\over {\ell}}\right]}(\beta) - {1\over 2}S_{\left[{{2h}\over {\ell}}\right]}(\beta)\right) - \left( 2\left[ {h\over {\ell}}\right] - {1\over 2} \left[ {{2h}\over {\ell}}\right] \right) - 2C_{\left[{h\over {\ell}}\right]}(\beta) + {1\over 2} C_{\left[{{2h}\over {\ell}}\right]}(\beta) + 
$$
$$
+ \left( 4\left\{ {h\over {\ell}}\right\}C_{\left[{h\over {\ell}}\right]}(\beta) - \left\{ {{2h}\over {\ell}}\right\}C_{\left[{{2h}\over {\ell}}\right]}(\beta)\right), 
$$
\par
i.e. 
$$
\enspace 
\Sigma_{\ell} = {{2\sin^2 \pi \beta \left[{h\over {\ell}}\right] - {1\over 2}\sin^2 \pi \beta \left[{{2h}\over {\ell}}\right]}\over {\sin^2(\pi \beta)}}\ell - h + \left( 2\left\{ {h\over {\ell}}\right\} \left( 1 + 2C_{\left[{h\over {\ell}}\right]}(\beta)\right) - {1\over 2} \left\{ {{2h}\over {\ell}}\right\} \left( 1 + 2C_{\left[{{2h}\over {\ell}}\right]}(\beta)\right) \right)\ell;
$$
\par
hence,
$$
\thinspace \enspace \thinspace 
{1\over {\ell }}\sum_{a}W(a\ell)e(a\beta) = {{4\sin^2 \pi \beta \left[{h\over {\ell}}\right] - \sin^2 \pi \beta \left[{{2h}\over {\ell}}\right]}\over {\sin^2(\pi \beta)}} + \left( 4\left\{ {h\over {\ell}}\right\} \sum_{0\le |a|\le {h\over {\ell}}}e(a\beta) - \left\{ {{2h}\over {\ell}}\right\} \sum_{0\le |a|\le {{2h}\over {\ell}}}e(a\beta) \right);
$$
\par
we distinguish two cases: first, $0\le \left\{ {h\over {\ell}}\right\} < {1\over 2}$ and, then, ${1\over 2}\le \left\{ {h\over {\ell}}\right\} < 1$. In the first, we have $\left[{{2h}\over {\ell}}\right]=2\left[{h\over {\ell}}\right]$:
$$
{1\over {\ell}} \sum_{a}W(a\ell)e(a\beta) = {{4\sin^4 \pi \beta \left[{h\over {\ell}}\right]}\over {\sin^2(\pi \beta)}} + \left( 4\left\{ {h\over {\ell}}\right\} \sum_{0\le |a|\le {h\over {\ell}}}e(a\beta) - \left\{ {{2h}\over {\ell}}\right\} \sum_{0\le |a|\le {{2h}\over {\ell}}}e(a\beta) \right), 
$$
\par
while in the second case we have $\left[{{2h}\over {\ell}}\right]=2\left[{h\over {\ell}}\right]+1$, so join (only for $2\left\{ {h\over {\ell}}\right\} - \left\{ {{2h}\over {\ell}}\right\}=1$) the term 
$$
-\cos 4\pi \beta \left[ {h\over {\ell}}\right] - \sin 4\pi \beta \left[ {h\over {\ell}}\right] \cot(\pi \beta) = 
\left( \hbox{\rm use\thinspace the\thinspace formula\thinspace for}\enspace C_{2\left[ {h\over {\ell}}\right]}(\beta),\;\hbox{\rm here} \right) 
$$
$$
= -2\left( \sum_{a\le 2\left[ {h\over {\ell}}\right]}\cos(2\pi a\beta) + {1\over 2}\right) 
= -\sum_{0\le |a|\le 2\left[ {h\over {\ell}}\right]}e(a\beta) 
= -\left( 2\left\{ {h\over {\ell}}\right\} - \left\{ {{2h}\over {\ell}}\right\} \right) \sum_{0\le |a|\le 2\left[ {h\over {\ell}}\right]}e(a\beta). \enspace \square
$$

\vfil
\eject

\par
\noindent {\bf 3. Proof of the Theorem.}
\par
\noindent
We will ignore the \thinspace $R_g(N,h)$ \thinspace that are \thinspace $\EssBdd Nh+h^3$(Good remainders!). Linking the Lemmas, 
$$
I_f(N,h) = \sum_{\ell \le 2h}\doublesum_{(d,q)=1}g(\ell d)g(\ell q){1\over q}\sum_{j\neq 0}\sum_{m\sim {N\over {\ell d}}}\cos {{2\pi jdm}\over q}\sum_{a\neq 0}W(a\ell)e_q(ja) 
$$
\par
\noindent
(save \thinspace $\EssBdd R_g(N,h)$, hereon); and using Lemma 2 instead of Lemma 1, see the introduction, 
$$
J_f(N,h) = \sum_{\ell \le 2h}\doublesum_{(d,q)=1}g(\ell d)g(\ell q){1\over q}\sum_{j\neq 0}\sum_{m\sim {N\over {\ell d}}}\cos {{2\pi jdm}\over q}\sum_{a\neq 0}S(a\ell)e_q(ja) 
$$
\par
\noindent 
Then, we'll show, for each $K$ like in Lemma 3, supported in $[-2h,2h]$, where uniformly bounded as $K\ll h$, 
$$
T_g(N,h)\defineq \sum_{\ell \le 2h}\doublesum_{(d,q)=1}g(\ell d)g(\ell q){1\over q}\sum_{j\neq 0}\sum_{m\sim {N\over {\ell d}}}\cos {{2\pi jdm}\over q} \sum_{a}K(a\ell)e_q(ja) \EssBdd Qh^2 + Q^2 h + R_g(N,h)
\leqno{(*)}
$$
\par
\noindent
In fact, we reintroduce terms with \thinspace $a=0$ (here \thinspace $K(0)=2h$), with contributes (${\sum}^{\ast}=$coprime to $d$)
$$
2h\sum_{\ell \le 2h}\doublesum_{(d,q)=1}g(\ell d){{g(\ell q)}\over q}\sum_{j\neq 0}\sum_{m\sim {N\over {\ell d}}}e_q(jdm) 
= 2h\sum_{\ell \le 2h}\sum_{d}g(\ell d)\sum_{m\sim {N\over {\ell d}}}\left({\sum_{q|dm}}^{*}g(\ell q) - {\sum_{q}}^{*}{{g(\ell q)}\over q}\right) 
\EssBdd Nh.
$$
\par
\noindent
(Once more from orthogonality, see above)\enspace We'll prove now (0). We may also join $j=0$ whenever $K=W$ : 
$$
\sum_{\ell \le 2h}\doublesum_{(d,q)=1}g(\ell d)g(\ell q){1\over q}\sum_{m\sim {N\over {\ell d}}} \sum_{a}W(a\ell) 
\EssBdd \sum_{\ell \le 2h}\doublesum_{(d,q)=1}{1\over q}\left({N\over {\ell d}}+1\right)h 
\EssBdd Nh
$$
\par
\noindent
and using (2), see Lemma 4, we get (only for $K=W$)
$$
T_g(N,h)=2\sum_{\ell \le 2h}\ell \doublesum_{(d,q)=1}g(\ell d)g(\ell q){1\over q}\sum_{j}\sum_{m\sim {N\over {\ell d}}}\cos {{2\pi jdm}\over q} 
\left( 1-\cos{{{2\pi j}\over q}[{h\over {\ell}}]}\right)\sum_{0\le |a|\le {h\over {\ell}}}\left( \left[ {h\over {\ell}}\right]-|a|\right)e_q(ja) + 
$$
$$
+ \sum_{\ell \le 2h}\ell \doublesum_{(d,q)=1}g(\ell d)g(\ell q){1\over q}\sum_{j}\sum_{m\sim {N\over {\ell d}}}\cos {{2\pi jdm}\over q} B\left( {h\over {\ell}}\right)\sum_{0\le |a|\le {{2h}\over {\ell}}}U_a\left( {h\over {\ell}}\right)e_q(-ja), 
$$
\par
\noindent
(plus negligible remainders), where $B$ and $U_a$ (uniformly on $a$) are bounded functions. From orthogonality, 
$$
{1\over q}\sum_{j}\sum_{m\sim {N\over {\ell d}}}\cos {{2\pi jdm}\over q} B\left( {h\over {\ell}}\right)\sum_{0\le |a|\le {{2h}\over {\ell}}}U_a\left( {h\over {\ell}}\right)e_q(-ja) = B\left( {h\over {\ell}}\right)\sum_{0\le |a|\le {{2h}\over {\ell}}}U_a\left( {h\over {\ell}}\right)\sum_{{m\sim {N\over {\ell d}}}\atop {m\equiv a(q)}}1, 
$$
\par
\noindent
whence we get (0), applying orthogonality (and the Lemmas) also on the main term, since for remainders we obtain: 
$$
\sum_{\ell \le 2h}\ell \doublesum_{(d,q)=1}g(\ell d)g(\ell q){1\over q}\sum_{j}\sum_{m\sim {N\over {\ell d}}}\cos {{2\pi jdm}\over q} B\left( {h\over {\ell}}\right)\sum_{0\le |a|\le {{2h}\over {\ell}}}U_a\left( {h\over {\ell}}\right)e_q(-ja)
\EssBdd Nh.
$$
\par
\noindent
We pass to the two bounds for our integrals. Now, from (1) and the well-known formula (Fej\'{e}r kernel) 
$$
\sum_{a}S(a)e_q(ja) 
= \sum_{0\le |a|\le 2h}(2h-|a|)e_q(ja) 
= {{\sin^2{{2\pi jh}\over q}}\over {\sin^2{{\pi j}\over q}}}, \enspace \hbox{\rm which\thinspace gives}\enspace \sum_{a}S(a\ell)e_q(ja)\ge 0 \qquad \forall j\neq 0
$$
\par
\noindent
(like in Lemma 4 proof), we have $\forall j\neq 0$, say, (for both $K=W,S$)
$$
\widehat{K}\left({j\over q}\right)\defineq \sum_{a}K(a\ell)e_q(ja)\ge 0 
$$
\par
\noindent
($\widehat{W}(0)\ge 0$ and $\widehat{S}(0)={{4h^2}\over {\ell}}+{\cal O}(h)$, trivially); whence (apart from $\EssBdd R_g(N,h)$), writing \lq \lq $\ast$\rq \rq \thinspace  for $(d,q)=1$,  
$$
T_g(N,h) = \sum_{\ell \le 2h}\doublesum_{(d,q)=1}g(\ell d)g(\ell q){1\over q}\sum_{j\neq 0}\sum_{m\sim {N\over {\ell d}}}\cos {{2\pi jdm}\over q} \widehat{K}\left({j\over q}\right) 
\EssBdd \max_{D\le Q} \sum_{\ell \le 2h} \sum_{q\sim {D\over {\ell}}}{1\over q} {\sum_{d\le 2q}}^{*} \sum_{j\neq 0} {1\over {\left \Vert {{jd}\over q}\right \Vert}} \widehat{K}\left({j\over q}\right) 
$$
\par
\noindent
due to [D, ch.\thinspace 25] 
$$
\sum_{m\sim {N\over {\ell d}}}e_q(jdm)\ll {1\over {\left \Vert {{jd}\over q}\right \Vert}}, 
$$
\par
\noindent
having used a dissection argument (over both $d,q$). Changing variables (with $\overline{d}d\equiv 1(q)$, here) into 
$$
\sum_{j\neq 0} {1\over {\left \Vert {{jd}\over q}\right \Vert}} \widehat{K}\left({j\over q}\right) 
= \sum_{0<|j|\le {q\over 2}} {1\over {\left \Vert {j\over q}\right \Vert}} \widehat{K}\left({{j\overline{d}}\over q}\right) 
\buildrel{K\,\hbox{\stampatello even}}\over{=\!=\!=\!=}
2q\sum_{j\le q/2}{1\over j}\widehat{K}\left({{j\overline{d}}\over q}\right) 
$$
\par
\noindent
gives 
$$
T_g(N,h)\EssBdd \max_{D\le Q} \sum_{\ell \le 2h} \sum_{q\sim {D\over {\ell}}}\sum_{j\le q/2}{1\over j}{\sum_{n\le 2q}}^{*} \widehat{K}\left({jn\over q}\right) 
$$
\par
\noindent
which, since 
$$
{\sum_{n\le 2q}}^{*} \widehat{K}\left({jn\over q}\right) \ll \sum_{n\le 2q}\widehat{K}\left({jn\over q}\right) 
= \sum_{a}K(a\ell)\sum_{n\le 2q}e_q(jan)=2q\sum_{{a}\atop {ja\equiv 0(q)}}K(a\ell) 
$$
\par
\noindent
and \lq \lq flipping\rq \rq \thinspace the divisors (say, {\stampatello f}, i.e., change $t$ into $q/t$) in the following 
$$
2q\sum_{j\le q/2}{1\over j}\sum_{{a}\atop {ja\equiv 0(q)}}K(a\ell) 
= 2\sum_{{t|q}\atop {t<q}}q\sum_{{j\le q/2}\atop {(j,q)=t}}{1\over j}\sum_{a\equiv 0(q/t)}K(a\ell) 
\buildrel{\hbox{\stampatello f}}\over{=\!=}
2\sum_{{t|q}\atop {t>1}}t\sum_{{j\le t/2}\atop {(j,t)=1}}{1\over j}\sum_{a\equiv 0(t)}K(a\ell) =
$$
$$
= 4h\sum_{{t|q}\atop {t>1}}t\sum_{{j\le t/2}\atop {(j,t)=1}}{1\over j} 
+2\sum_{{t|q}\atop {1<t\le 2h}}t\sum_{{j\le t/2}\atop {(j,t)=1}}{1\over j}\sum_{{a\neq 0}\atop {a\equiv 0(t)}}K(a\ell) 
\EssBdd qh+h^2, 
$$
\par
\noindent
finally entails 
$$
T_g(N,h) \EssBdd Nh + \max_{D\le Q} \sum_{\ell \le 2h} \sum_{q\sim {D\over {\ell}}}(qh+h^2) 
\EssBdd Nh + Q^2h + Qh^2 \Rightarrow (\ast).\enspace \square 
$$
\medskip
\par
\noindent
\centerline{\bf References}
\smallskip
\item{\bf [C-S]} Coppola, G. and Salerno, S.\thinspace - \thinspace {\sl On the symmetry of the divisor function in almost all short intervals} \thinspace - \thinspace Acta Arith. {\bf 113} (2004), {\bf no.2}, 189--201. $\underline{\tt MR\enspace 2005a\!:\!11144}$
\smallskip
\item{\bf [D]} \thinspace Davenport, H.\thinspace - \thinspace {\sl Multiplicative Number Theory} \thinspace - \thinspace Third Edition, GTM 74, Springer, New York, 2000. $\underline{{\tt MR\enspace 2001f\!:\!11001}}$
\smallskip
\item{\bf [L]} \thinspace Linnik, Ju.V.\thinspace - \thinspace {\sl The Dispersion Method in Binary Additive Problems} \thinspace - \thinspace Translated by S. Schuur \thinspace - \thinspace American Mathematical Society, Providence, R.I. 1963. $\underline{\tt MR\enspace 29\# 5804}$
\smallskip
\item{\bf [T]} \thinspace Tenenbaum, G.\thinspace - \thinspace {\sl Introduction to Analytic and Probabilistic Number Theory} \thinspace - \thinspace Cambridge Studies in Advanced Mathematics, {\bf 46}, Cambridge University Press, 1995. $\underline{\tt MR\enspace 97e\!:\!11005b}$
\smallskip
\item{\bf [V]} \thinspace Vinogradov, I.M.\thinspace - \thinspace {\sl The Method of Trigonometrical Sums in the Theory of Numbers} - Interscience Publishers LTD, London, 1954. $\underline{{\tt MR \enspace 15,941b}}$
\medskip
\leftline{\tt Dr.Giovanni Coppola}
\leftline{\tt DIIMA - Universit\`a degli Studi di Salerno}
\leftline{\tt 84084 Fisciano (SA) - ITALY}
\leftline{\tt e-mail : gcoppola@diima.unisa.it}

\bye